\documentclass[oneside,reqno,12pt]{amsart}
\usepackage[greek,english]{babel}
\usepackage{amsthm}
\usepackage{amsbsy}
\usepackage{amsfonts}
\usepackage{graphicx}
\usepackage{hyperref}
\hypersetup{colorlinks=true, citecolor=blue, linkcolor=red}

\newtheorem{thm}{Theorem}
\newtheorem{cor}{Corollary}

\newtheorem{rem}{Remark}

\numberwithin{equation}{section} \numberwithin{lem}{section}
\numberwithin{thm}{section} \numberwithin{cor}{section}
\numberwithin{pro}{section} \numberwithin{rem}{section}

\begin{document}

\title[Instability and nonordering of localized steady states]{Instability and nonordering of localized steady states to a classs of reaction-diffusion equations in $\mathbb{R}^N$}
\author{Christos Sourdis}
\address{National and Kapodistrian University of Athens, Department of Mathematics, Athens, Greece.
}
\email{sourdis@uoc.gr}

\date{\today}
\begin{abstract}
 We show that the elliptic problem $\Delta u+f(u)=0$ in $\mathbb{R}^N$, $N\geq 1$, with $f\in C^1(\mathbb{R})$   and $f(0)=0$ does not have nontrivial stable solutions that decay to zero at infinity, provided that $f$ is nonincreasing near the origin. As a corollary, we can show  that any two nontrivial solutions that decay to zero at infinity must intersect each other, provided that at least one of them is sign-changing.
This property was previously known only in the case where both solutions are positive with a different approach.
We also discuss implications of our main result on the existence of monotone heteroclinic solutions to the corresponding reaction-diffusion equation.
 \end{abstract}
 \maketitle

\section{Introduction}\subsection{The setting and known results}
We consider the elliptic problem
 \begin{equation}\label{eqsteady}
   \Delta u+f(u)=0\ \textrm{in}\ \mathbb{R}^N,\ N\geq 1,
 \end{equation}
  with \begin{equation}\label{eqf}f\in C^1(\mathbb{R})\ \textrm{and}\ f(0)=0.\end{equation} Solutions will be understood in the classical sense (this property is valid for  locally bounded, distributional solutions   by standard elliptic regularity theory). A solution is called \emph{localized} if it satisfies
  \begin{equation}\label{eqlocal}
    u(x)\to 0\ \textrm{as}\ |x|\to \infty.
  \end{equation} The above problem  has been studied extensively for various nonlinearities (see for example \cite{hamel} and the references therein).
A solution is called \emph{stable} if
\begin{equation}\label{eqviol}
  \int_{\mathbb{R}^N}^{}\left\{|\nabla \varphi|^2-f'(u)\varphi^2 \right\}dx\geq 0\ \ \forall\ \varphi\in C_0^\infty(\mathbb{R}^N),
\end{equation}
see for instance \cite{dupaigne}. Otherwise, it is called \emph{unstable}.

Let us briefly mention some important rigidity results for stable solutions to (\ref{eqsteady}) that hold without any further restrictions on $f$ and without the assumption that $f(0)=0$. As was shown in \cite{cabre}, there are no nontrivial stable solutions   in the Sobolev space $W^{1,2}(\mathbb{R}^N)$, $N\geq 1$. In fact, it was observed therein that if $N\geq 2$ the same conclusion holds
under the weaker assumption that $|\nabla u|\in L^2(\mathbb{R}^N)$ (under the further assumption that $u$ has finite energy,
this  was proven much earlier in \cite{derrick}). If $N=2$, then bounded stable solutions depend only on one variable (possibly after a rotation), see \cite{farina,dancer}.  If $N\leq 10$, there are no nonconstant radially symmetric stable solutions that are bounded    (see \cite{cabre,vvv}).
If one further assumes that $f\geq 0$  and $N\leq 10$, then there are no nonconstant stable solutions  that are bounded from below, see \cite{dp} which relies on the a-priori estimates that were obtained in \cite{figalli}.

Positive solutions of (\ref{eqsteady})-(\ref{eqlocal}), under the assumption (\ref{eqf}) on $f$, are known to be radially symmetric and decreasing with respect to some point, provided that
\begin{equation}\label{eqf'}
  f'(s)\leq 0\ \textrm{for small}\ |s|,
\end{equation}
(see \cite{ln}). Moreover, in this case any two such positive solutions must intersect  (see \cite[Lem. 3.2]{busca}).

If
\begin{equation}\label{eqstrauss}
f'(0)<0,
\end{equation}
then any localized solution to (\ref{eqsteady}) and its gradient must decay  exponentially fast as $|x|\to \infty$ (by a standard barrier argument). Thus, by the above discussion, such  solutions must be unstable (see also \cite{strauss} for the case of radial solutions).
In this paper, we will show that this property continues to hold under the weaker condition (\ref{eqf'}).
We note in passing that it was shown in \cite{villegas2020}, under the sole assumption (\ref{eqf}), that there are no nontrivial localized solutions that are minimizers in the sense of Morse (a property that is stronger than stability). Moreover, as a corollary of our main result, we will generalize the aforementioned result of \cite{busca}.

\subsection{Our results and methods of proof}

Our main result is the following.
\begin{thm}\label{thm0}
If $u$ is a stable solution to (\ref{eqsteady}) and (\ref{eqlocal}) with $f$ satisfying (\ref{eqf}) and (\ref{eqf'}), then $u\equiv 0$.
\end{thm}

In contrast, if $f(u)=|u|^{p-1}u$, $p>1$, and $N\geq 11$, it was shown in \cite{guiComptia} by Hardy's inequality that the  positive   radial solutions of (\ref{eqsteady}) and (\ref{eqlocal}) are stable provided that
$p$ is greater than or equal to  the Joseph-Lundgren  exponent. It is worth mentioning that their stability also follows from the fact that they form a foliation
and the convexity of the nonlinearity as in \cite[Prop. 1.3.2]{dupaigne} (see also \cite{guiComptia2} for this viewpoint). The point is that such pure power nonlinearities satisfy (\ref{eqf}) but not (\ref{eqf'}) (we refer to \cite{cabre} and \cite{vvv} for related examples).

Our proof
proceeds by showing that all first order partial derivatives of $u$ are identically equal to zero. These satisfy the linearized equation  on $u$ and tend to zero at infinity (by standard elliptic estimates).
 It is a well known fact that the stability of $u$ implies the existence of a positive solution $\Psi$ to the aforementioned linearized equation. If $\liminf_{|x|\to \infty}\Psi>0$, then
 the linearized operator satisfies the maximum principle (in the whole space, see \cite{bcv}), and the assertion of the theorem follows at once. In any case, the existence of $\Psi>0$ implies that
 the maximum principle holds in any bounded domain. On the other hand, the assumption (\ref{eqf'}) implies that the maximum principle holds in the exterior of large balls around the origin. Remarkably, these two separate properties can be combined to show that the linearized operator satisfies the maximum principle in the whole space, and therefore one can conclude as before. This can be shown by adapting an argument from \cite{cabreU} which is based on considering the quotient of the solution over $\Psi$. However, we found it more convenient to argue directly using Serrin's sweeping principle, in the spirit of the moving plane argument of \cite{ln}.

As we have already mentioned, for $f$  as in Theorem \ref{thm0},  it was shown in \cite[Lem. 3.2]{busca} that any two positive solutions of (\ref{eqsteady}) and (\ref{eqlocal})
must intersect each other. This was accomplished by the famous sliding method \cite{bn}, exploiting that (\ref{eqsteady}) is invariant under translations. However, this approach breaks down in the case of sign-changing solutions. On the other hand, our Theorem \ref{thm0} can be used to show that    an analogous intersection property   holds in this case as well.
More precisely, the following corollary complements the aforementioned result of \cite{busca} for positive solutions (our proof also works in this case but is not as direct).

\begin{cor}\label{cor0}Let $f$ satisfy (\ref{eqf}) and (\ref{eqf'}).
If $u_1$ and $u_2$ are two distinct nontrivial solutions of  (\ref{eqsteady}) and (\ref{eqlocal}) with at least one of them being sign-changing, then $u_1-u_2$ must change sign.
\end{cor}

The main idea of the proof is the following. If they were ordered,  since both are unstable by Theorem \ref{thm0}, we can use a dynamical systems  argument
to show that there exists a stable solution between them (note that it must be nontrivial), which is in contradiction to Theorem \ref{thm0}. This type of arguments are well known in the case of bounded domains
(see \cite{pao,sat}). The case of the whole space requires a bit of extra care.

\subsection{Applications}
Let us briefly discuss some interesting implications of Theorem \ref{thm0}.

Let $f\in C^1(0,\infty)\cap C[0,\infty)$ satisfy $f(0)=0$ and (\ref{eqf'}). Assume that there exists a positive solution $w$ to (\ref{eqsteady}) and (\ref{eqlocal}).
As we have already mentioned, $w$ has to be radially symmetric and decreasing with respect to some point and there does not exist another such solution below it.
Moreover, we know from Theorem \ref{thm0} that $w$ is unstable. So, there exists a large ball  $B_R$ such that the principal eigenvalue of
\begin{equation}\label{eqEVs}
-\Delta \psi-f'(w)\psi=\lambda \psi \ \textrm{in}\ B_R;\ \psi=0 \ \textrm{on}\ \partial B_R,
\end{equation}
is negative (one can take $B_R$ to include the support of a test function that violates (\ref{eqviol})).
Actually, in  \cite{busca}, \cite[Sec. 1.4]{hamel} the existence of such a ball was shown under the stronger condition (\ref{eqstrauss}) and an approximation argument. Armed with the above information in $B_R$,  the approach in the latter reference applies to establish that the reaction-diffusion equation
\[
u_t=\Delta u+f(u), \ x\in \mathbb{R}^N, \ t\in \mathbb{R},
\]
admits a heteroclinic solution such that $u_t<0$ and
 \[
  u\to w\ \textrm{as}\ t\to -\infty;\ u\to 0\ \textrm{as}\ t\to +\infty,\ \textrm{uniformly in}\ \mathbb{R}^N.
 \]
\begin{rem}
 It is easy to see that in  the     scheme of \cite{hamel} for the construction of the aforementioned heteroclinic solution one can also take as  initial condition the function
 \begin{equation}\label{eqpsar}
    \min\{w(\cdot+\varepsilon e),w\}, \ \textrm{where}\ e=(1,0,\cdots,0)\ \textrm{and}\ 0<|\varepsilon|\ll 1,
 \end{equation}
(from the aforementioned result of \cite{busca},  $w$ must intersect any translate of itself),
 instead of
 \begin{equation}\label{eqsupersol}
w_\varepsilon := \left\{
 \begin{array}{ll}
    w-\varepsilon \Phi, & x\in B_R, \\
   w, & x\in \mathbb{R}^N\setminus B_R,
 \end{array}
 \right.
 \end{equation} where  $\Phi>0$  stands for the principal eigenfunction of (\ref{eqEVs}) and $0<\varepsilon \ll 1$,
 which was used therein. The point is that both are strict weak supersolutions to (\ref{eqsteady}) (the positivity of $w$ is not used here).

 This observation, which implies that $w$ is a dynamically unstable steady state in $L^\infty(\mathbb{R}^N)$ (without any a-priori information on the linearized operator), was actually our heuristic motivation behind
 Theorem \ref{thm0}. We point out that the solution of the corresponding Cauchy problem with initial condition as in (\ref{eqpsar}) or (\ref{eqsupersol}) converges monotonically to zero as $t\to +\infty$ (see \cite{busca,hamel}).
 \end{rem}
\subsection{Outline of the paper} The rest of the paper is devoted to the proofs of Theorem \ref{thm0} and Corollary \ref{cor0}.
\section{Proofs}\subsection{Proof of Theorem \ref{thm0}}
\begin{proof}
Our goal is to prove that \begin{equation}\label{eqGoal}
                            \partial_{x_i} u\equiv0\  \textrm{for}\  i=1,\cdots,N,
                          \end{equation}from where the assertion of the theorem follows at once. We first note that standard elliptic estimates yield \begin{equation}\label{eqyela}
                                                                                                \partial_{x_i} u\to 0\ \textrm{as}\ |x|\to \infty \  \textrm{for}\  i=1,\cdots,N.
                                                                                              \end{equation}
Moreover, each $\partial_{x_i}u\in C^{1,\alpha}(\mathbb{R}^N)$, $\alpha\in (0,1)$, satisfies weakly the linearized equation of (\ref{eqsteady}) at $u$.

Since $u$ is stable, as in \cite[Prop. 4.2]{alberti} or \cite[Thm. 1.7]{bcnp}, there exists a $\Psi \in C^{1,\alpha}(\mathbb{R}^N)$, $\alpha\in (0,1)$, such that
\[
-\Delta \Psi-f'(u)\Psi=0\ \textrm{(in the weak sense) and}\ \Psi>0\ \textrm{in}\ \mathbb{R}^N.
\]
We will show that (\ref{eqGoal}) holds with the use of Serrin's sweeping principle  (see \cite[Thm. 2.7.1]{sat}).
For a fixed $i\in \{1,\cdots,N\}$, let us consider the set
\[\Lambda = \left\{
\lambda \geq 0 \  : \
\partial_{x_i} u \leq \mu \Psi \ \textrm{in}\ \mathbb{R}^{N}
\ \textrm{for every}\ \mu \geq \lambda
\right\}.
\]	
Our goal is to show that $\Lambda = [0, \infty)$, which will yield $\partial_{x_i} u \leq 0$. We can
also apply the same argument, with $\partial_{x_i} u$ replaced by $-\partial_{x_i} u$, to obtain $\partial_{x_i} u \geq 0 $ and therefore
conclude.

By virtue of (\ref{eqlocal}) and (\ref{eqf'}), there exists an $R>0$ such that
\begin{equation}\label{eqf''}
f'(u)\leq 0\ \textrm{for}\ |x|\geq R.
\end{equation}
Clearly, there exists a $\bar{\lambda}>0$ such that
\[
\partial_{x_i} u \leq \bar{\lambda} \Psi \ \textrm{for}\ |x|\leq R.
\]
Since both $\partial_{x_i}u$ and $\Psi$ satisfy the  linearized equation of (\ref{eqsteady}) at $u$, it follows from  (\ref{eqyela}), (\ref{eqf''}) and the maximum principle (for weak solutions, see \cite[p. 48]{bcv})  that the above ordering is also valid for $|x|>R$, i.e. $\bar{\lambda}\in \Lambda$.
Hence, $\Lambda$ is an interval of the form $[\tilde{\lambda},\infty)$ for some $\tilde{\lambda}\in [0,\bar{\lambda}]$.

It remains to show that $\tilde{\lambda}=0$. To this end, we will argue by contradiction and suppose that $\tilde{\lambda}>0$.
From the relation
\[
\partial_{x_i} u \leq \tilde{\lambda} \Psi \ \textrm{in}\ \mathbb{R}^N,
\]
and the strong maximum principle which implies that the above inequality is strict (the possibility that $\partial_{x_i}u\equiv \tilde{\lambda} \Psi$ is easily excluded from (\ref{eqlocal})), we infer that there exists a $\delta\in (0,\tilde{\lambda}/2)$ such that
\[
\partial_{x_i} u \leq (\tilde{\lambda}-\delta) \Psi \ \textrm{for}\ |x|\leq R.
\]
Then, as before, we deduce by the maximum principle that the above relation holds in $\mathbb{R}^N$, which  contradicts the minimality of $\tilde{\lambda}$ and completes the proof of the theorem.
\end{proof}

\subsection{Proof of Corollary \ref{cor0}}
\begin{proof} Without loss of generality, we may assume  that $u_2$ is sign-changing.
Let us argue  by contradiction and suppose that the assertion of the corollary is false, i.e. $u_1$ and $u_2$ are ordered. So, thanks to the strong maximum principle, we get $u_1<u_2$ in $\mathbb{R}^N$ (again without loss of generality). Since Theorem \ref{thm0} guarantees that
$u_2$ is unstable, as we have already mentioned, there exists a strict weak supersolution $u_{2,\varepsilon}$  to (\ref{eqsteady}) of the form (\ref{eqsupersol}) such that $u_1<u_{2,\varepsilon}\leq u_2$. In fact, by taking a smaller $\varepsilon>0$ and a larger radius $R>0$ if necessary, we have that $u_{2,\varepsilon}$ is sign-changing.
Let us now consider the solution $v$ of the Cauchy problem
\[
\left\{\begin{array}{ll}
         u_t=\Delta u+f(u), & x\in \mathbb{R}^N,\ t>0, \\
         u(x,0)=u_{2,\varepsilon}(x), &x\in \mathbb{R}^N.
       \end{array}
 \right.
\]
It is well known that $v_t<0$ (see for instance \cite[Prop. 52.20]{qs}), and thus
\[
v(\cdot,t)\to z(\cdot) \ \textrm{as}\ t\to +\infty,\ \textrm{uniformly in}\ \mathbb{R}^N,
\]
where $z$ is a steady state such that $u_1\leq z<u_{2,\varepsilon}$ in $\mathbb{R}^N$.
Similarly, taking into account that $z$ cannot be identically equal to zero (recall that $u_{2,\varepsilon}$ is sign-changing),  there exists a strict weak lower solution $z_{\epsilon}$ of the form (\ref{eqsupersol}) with $0<-\epsilon\ll 1$ such that \begin{equation}\label{eqstrict2}
                                                                                                                z\leq z_{\epsilon}<u_{2,\varepsilon}\ \ \textrm{in}\ \mathbb{R}^N.
                                                                                                              \end{equation}

On the other hand, the maximum principle yields
\[
z_\epsilon (x)<v(x,t),\ x\in \mathbb{R}^N,\ t>0.
\]
Letting $t\to +\infty$ in the above relation gives
\[
z_\epsilon\leq z \ \ \textrm{in}\ \mathbb{R}^N,
\]
which together with (\ref{eqstrict2}) contradicts the fact that $z_\epsilon$ is a strict lower solution.
\end{proof}




\subsection*{Acknowledgments} The author would like to thank Prof. Hamel for bringing to his attention \cite{hamel}. Furthermore, he wishes to thank  the referee for offering pertinent remarks and suggestions. Moreover, he would like to thank IACM of FORTH, where
this paper was written, for the hospitality. This work has received funding from
the Hellenic Foundation for Research and Innovation (HFRI) and the General
Secretariat for Research and Technology (GSRT), under grant agreement No
1889.

\end{document}